\documentclass[12pt]{article}
\usepackage[utf8]{inputenc}
\usepackage{amsmath, amssymb, amsthm}
\usepackage{graphicx}
\usepackage{color}
\usepackage{subfigure}
\usepackage{lineno}
\usepackage{amsmath,amssymb}

\usepackage{amsthm}
\usepackage{amsfonts}
\usepackage{mathrsfs}
\usepackage{indentfirst}
\usepackage{eucal}
\usepackage{eufrak}
\usepackage{color, graphicx}
\usepackage{listings}
\usepackage{times}
\usepackage[english]{babel}
\usepackage{fancyhdr}
\usepackage{makeidx}
\newtheorem{theorem}{Theorem}[section]
\newtheorem{corollary}{Corollary}

\newtheorem{lemma}[theorem]{Lemma}

\theoremstyle{definition}

\date{}
\usepackage[table,xcdraw]{xcolor}  
\usepackage{graphicx}  
\usepackage{mathtools}
\usepackage{tikz}
\begin{document}
\title{On connected graphs with finite spectral redundancy index and Pythagorean triplets}
\author{ Pawan Kumar$^{a}$, S. Pirzada$^{b}$, S. Merajuddin$^{a}$ \\~\\
{\em $^{a}$Department of Applied Mathematics, Aligarh Muslim University, Aligarh, India}\\
 \texttt{paone101@gmail.com},~ \texttt{meraj1957@rediffmail.com}\\
{\em $^{b}$Department of Mathematics, University of Kashmir, Srinagar, Kashmir, India}\\
\texttt{pirzadasd@kashmiruniversity.ac.in}}
\date{}

\pagestyle{myheadings} \markboth{Pirzada}{On connected graphs with finite spectral redundancy index}

\maketitle
\begin{abstract}
This article investigates spectral redundancy, a concept initially introduced by Alberto Seeger. Spectral redundancy arises when different connected induced subgraphs of a graph share the same spectral
radius in their adjacency spectrum. Let \(b(G)\) denote the total number of non-isomorphic induced subgraphs of \(G\), and \(c(G)\) represents the cardinality of the set of spectral radius of all connected induced subgraphs of \(G\). The spectral redundancy of a graph \( G \) is defined as the ratio \( \frac{b(G)}{c(G)} \). The supremum of this ratio across all graphs in a family is called the spectral redundancy index of that family. We focus on a family of graphs that exhibit spectral redundancy and we find out the spectral redundancy index of this family. Furthermore, we investigate the connection between the spectral redundancy of these graphs and the presence of Pythagorean triplets.
\end{abstract}
\noindent{ \textbf{MSC 2020:} 05C50, 15A42

\section{Introduction}
The spectrum of a graph, comprising the eigenvalues of its adjacency matrix, is a powerful tool for understanding the structural properties of graphs. Spectral graph theory has significantly contributed to understanding graph properties such as connectivity, symmetry, graph coloring, and graph isomorphism detection \cite{ineq,induced_lemma}. Among the spectral parameters, the spectral radius, which is the largest eigenvalue of the adjacency matrix, holds particular significance due to its implications in various applications\cite {induced_lemma}.

Recently, a novel spectral concept related to the adjacency matrix has emerged, known as the complementarity spectrum of a graph. This concept was first introduced by Seeger \cite{intr} for general square matrices. Since then, it has been studied extensively in various contexts for its distinctive properties and applications.

Formally, let \( A \) be a real matrix of order \( n \). A scalar \( \lambda \in \mathbb{R} \) is called a \textit{complementarity eigenvalue} of \( A \) if there exists a vector \( x \in \mathbb{R}^n \setminus \{0\} \) such that
\[
x \geq 0, \quad A x - \lambda x \geq 0 \quad \text{and} \quad \langle x, A x - \lambda x \rangle = 0,
\]
where \( x \geq 0 \) indicates that \( x \) is componentwise nonnegative.

Let \( G \) be a connected graph of order \( n \) with adjacency matrix \( A_G \). Fernandes et al. \cite{Fer} extended the concept of complementarity spectrum to graphs. They defined the complementarity spectrum of a graph as the complementarity spectrum of its adjacency matrix \( A_G \).

An \emph{induced subgraph} of \( G \) is formed by selecting a subset of the vertices and including all the edges that connect pairs of vertices within this subset \cite{Alon}. Fernandes et al. \cite{Fer} also showed that the complementarity spectrum of \( G \), denoted by \( \Pi(G) \), is the set of Perron roots (spectral radius) of the adjacency matrices of all induced non-isomorphic connected subgraphs of \( G \). Mathematically, this is expressed as
\[
\Pi(G) = \{ \rho(F) : F \in \mathcal{S}(G) \},
\]
where \( \mathcal{S}(G) \) is the set of all induced connected subgraphs of \( G \) and \( \rho(F) \) denotes the spectral radius of \( F \). Unlike other spectrum of graph, the complementarity spectrum of a graph of order \( n \) always has more than \( n \) elements for non-elementary graphs (i.e., graphs other than paths, cycles, stars, and complete graphs).   Consequently, this larger range of possible complementarity eigenvalues makes the complementarity spectrum a powerful tool for analyzing graph structures and distinguishing between different graphs. Notably, no pair of non-isomorphic graphs of the same order has been found to share the same complementarity spectrum to date.

Let \( b(G) \) be the total number of non-isomorphic induced subgraphs of \( G \) and let \( c(G) \) denote the size of the complementarity spectrum \( \Pi(G) \). Fernandes et al. \cite{Fer} established the following inequalities
\[
c(G) \leq b(G) \quad \text{and} \quad n \leq b(G) \leq 2^n - 1.
\]
The lower bound for \( c(G)\geq n \) was later determined by Seeger \cite{Seeger}. In addition, Seeger made another important contribution \cite{See} by introducing the concept of the spectral redundancy of \( G \), defined as the ratio
\[
r(G) = \frac{b(G)}{c(G)}.
\]
Spectral redundancy \( r(G) \) measures the degree of repetition in spectral radius among all connected induced subgraphs of $G$. A high spectral redundancy indicates that a large number of induced subgraphs share the same spectral radius, whereas low redundancy suggests less repetition of spectral radius among induced subgraphs. For a possibly infinite family of connected graphs \( \mathbf{G} \), the spectral redundancy index of \( \mathbf{G} \) is defined by
\[ \mathfrak{r}(\mathbf{G}) = \underset{G \in \mathbf{G}}{\sup} \, r(G). \]
In his work \cite{See}, Seeger examined the spectral redundancy index for different families of graphs and demonstrated for instance that the set of all broom graphs has a spectral redundancy index of $4/3$. The foundational work by Seeger \cite{friendship}, Fernandes \cite{Fer2}, and others has established the basis for continued research in this dynamic area of mathematics. Subsequent studies by Seeger \cite{Seeger}, Seeger and Sossa \cite{see2}, \cite{Sossa} and Pinheiro et al. \cite{Pin} have further deepened and refined the understanding of this concept, contributing to the growth of research in the field.

In this article, we comprehensively examine the spectral redundancy within a special infinite family of connected graphs. We investigate the conditions under which two graphs in this family share a common spectral radius and explore how the spectral redundancy in this family related to the presence of Pythagorean triplets. A \textit{Pythagorean triplet} is defined as a set of three positive integers \((a, b, c)\) that satisfy the equation $a^2 + b^2 = c^2$, where \(a\) and \(b\) are referred to as the \textit{legs} and \(c\), the largest number, is termed the \textit{hypotenuse}. A Pythagorean triplet is said to be \textit{primitive} if \(a\), \(b\) and \(c\) are coprime integers. Throughout this article, we will use both the abbreviation PT and \textit{Pythagorean triplet} interchangeably for convenience.

The structure of the paper is as follows. Section 2 examines a particular family of connected graphs that exhibit repeated spectral radius among their members. Section 3 explores the spectral redundancy index, providing an in-depth discussion of its properties and implications in graph theory. Finally, Section 4 presents a detailed discussion that highlights the broader significance of the findings and potential directions for future research.

\section{Spectrally Redundant Graphs }

Let \( G_{p,q} \) denote a graph of order \( p+q+2 \), where a set of $p$ pendant vertices is connected to a central vertex, creating a star graph configuration around it. In addition, this central vertex is connected to another vertex through $q$ parallel paths, each consisting of two edges. For example, the graph \( G_{4,6} \) is depicted in Figure \ref{garl}. We define the set \( \mathcal{G}(p,q) \) as follows
\[
\mathcal{G}(p,q) = \{ G_{p,q} \mid p \geq 0, q \geq 1 \}.
\]
In this article, we refer to the members of the family $\mathcal{G}(p,q)$ as garlic graphs (having resemblance to a garlic-like structure) for clarity and ease of reference in our discussions and proofs.

A connected induced subgraph of \( G_{p,q} \) is either another garlic graph or a star graph. A star graph, \( S_r \), consists of one central vertex connected to \( r-1 \) pendant vertices. A star graph \( S_{r} \) is a special case of a \( G \)-graph with $p=r-1$ and \( q = 0 \), that is \( S_{r} \) $\cong$ \( G_{r-1,0} \). While this notation can be used interchangeably, it should be noted that \( G \)-graphs with \( q = 0 \) were not included in the original definition of the \( G \)-graph family as introduced. The set of all connected non-isomorphic induced subgraphs of a garlic graph \( G_{p,q} \) can be described by the following
\[ \mathcal{S}(G_{p,q}) = \{ G_{s,t} : 0 \leq s \leq p, 1 \leq t \leq q \} \cup \{ S_r : 1 \leq r \leq p+q+1 \}. \]
This expression indicates that for any garlic graph \( G_{p,q} \), the connected induced subgraphs include all possible garlic graphs \( G_{s,t} \) with \( s \) pendent vertices (where \( 0 \leq s \leq p \)) and \( t \) parallel paths (where \( 1 \leq t \leq q \)), as well as all possible star graphs with up to \( p+q+1 \) vertices.

We now determine the number of induced non-isomorphic connected subgraphs of a garlic graph, which will help in analyzing its spectral properties and overall graph behavior.
\begin{lemma}\label{subgraph}
    Let \( G_{p,q} \) be a garlic graph. The number of induced non-isomorphic connected subgraphs of \( G_{p,q} \), denoted by \( b(G_{p,q}) \), is given by
\[
b(G_{p,q}) = qp + 2q + p.
\]
\end{lemma}
\begin{proof}
    Consider the set of garlic graphs \( \{ G_{s,t} : 0 \leq s \leq p, 1 \leq t \leq q \} \), where \( s \) represents the number of pendent vertices and \( t \) represents the number of parallel paths. There are \( (p+1) \) choices for \( s \) and \( q \) choices for \( t \), giving a total of \( (p+1)q \) such subgraphs.

Next, consider the set of star subgraphs \( \{ S_r : 1 \leq r \leq p+q \} \). There are \( p+q+1 \) such subgraphs. Notably, the graph \( G_{0,1} \) is isomorphic to \( S_2 \), so we must subtract one to avoid double-counting. Therefore, the total number of non-isomorphic connected subgraphs of \( G_{p,q} \) is
\[
b(G_{p,q}) = q(p+1) + (p+q+1) - 1 = qp + 2q + p.
\]
This completes the proof.
\end{proof}
Having identified the non-isomorphic induced subgraphs of $G_{p,q}$, the next step is to analyze its spectral properties by deriving the characteristic equations of its adjacency matrix. The following lemma presents the characteristic equation of the adjacency matrix of graph \( G_{p,q} \).

\begin{figure}[h!]
\centering
\begin{tikzpicture}[scale=1.5, every node/.style={circle, draw, fill=black, inner sep=3pt}]
\coordinate (C) at (0, 0); 
\coordinate (L) at (0, -2); 
\coordinate (A1) at (-0.5, 1); 
\coordinate (A2) at (-1.25, 1); 
\coordinate (A3) at (0.5, 1);  
\coordinate (A4) at (1.25, 1);  
\coordinate (A5) at (2, 1);  
\coordinate (A6) at (3, 1);  
\coordinate (I1) at (-1.5, -1); 
\coordinate (I2) at (-1, -1);   
\coordinate (I3) at (-0.5, -1); 
\coordinate (I4) at (0.5, -1);  
\coordinate (I5) at (1, -1);    
\coordinate (I6) at (1.5, -1);  

\node[fill=black] (V) at (C) {};

\node[fill=black] (L) at (L) {};

\node[fill=black] (P1) at (A1) {};
\node[fill=black] (P2) at (A2) {};
\node[fill=black] (P3) at (A3) {};
\node[fill=black] (P4) at (A4) {};

\node[fill=black] (INT1) at (I1) {};
\node[fill=black] (INT2) at (I2) {};
\node[fill=black] (INT3) at (I3) {};
\node[fill=black] (INT4) at (I4) {};
\node[fill=black] (INT5) at (I5) {};
\node[fill=black] (INT6) at (I6) {};

\foreach \i in {1, 2, 3, 4} {
    \draw[line width=0.4 mm] (C) -- (P\i);
}

\foreach \i in {1, 2, 3, 4, 5, 6} {
    \draw[line width=0.4 mm] (C) -- (INT\i) -- (L);
}
\end{tikzpicture}
\caption{The garlic graph $G_{4,6}$ with 4 pendent vertices and 6 parallel path of length two}\label{garl}
\end{figure}
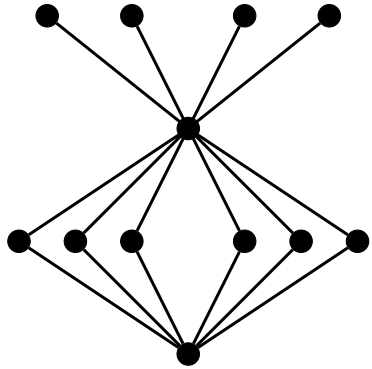
 \begin{lemma}
 The characteristic equation of the graph $G_{p,q}$ is given by the following equation
 \begin{equation}\label{char_eq}
 x^{p+q-2}\left[x^4-(p+2q)x^2+pq\right]=0.
 \end{equation}
 In particular, the spectral radius of \( G_{p,q} \) admits the explicit characterization
\begin{equation}\label{spectral_radius}
	\rho(G_{p,q}) = \sqrt{\frac{f(p, q)}{2}},
	\, \text{where}  \,
	f(p, q) = p + 2q + \sqrt{p^2 + 4q^2}.
\end{equation}
 \end{lemma}
\begin{proof}
	We observe that the graph \( G_{p,q} \) is constructed by coalescing two well-known graphs, the complete bipartite graph \( G_{0,q} = K_{2,q} \) and the star graph \( S_{p+1} \). Specifically, the central vertex of the star is merged with one of the two high-degree vertices of \( K_{2,q} \).
	
	A standard result on the characteristic polynomial of coalesced rooted graphs (see Corollary 2(b) in \cite{schwenk}) states that if two graphs \( F \) and \( H \) are coalesced at specific vertices, their characteristic polynomial satisfies the relation
	\[
	\phi(G_{p,q}, x) = \phi(F, x) \phi(H - u, x) + \phi(F - v, x) \phi(H, x) - x \phi(H - u, x) \phi(F - v, x),
	\]
	where \( u \) and \( v \) are the coalesced vertices, and \( \phi(H - u, x) \) and \( \phi(F - v, x) \) denote the characteristic polynomials of the respective graphs with those vertices removed.
	Using the known characteristic polynomials of \( K_{2,q} \) and \( S_{p+1} \), we obtain
	\[
	\phi(G_{p,q}, x) = x^{p+q-2} \left( x^4 - (p + 2q)x^2 + pq \right).
	\]
	This establishes the characteristic equation of \( G_{p,q} \).
	\vspace{0.2cm}
	
	The spectral radius is the largest eigenvalue and thus corresponds to the largest root of the biquadratic equation
	\[
	x^4 - (p+2q)x^2 + pq = 0.
	\]
	Solving for \( x \) yields the spectral radius
	\[
	\rho(G_{p,q}) = \left[\frac{p+2q + \sqrt{p^2 + 4q^2}}{2}\right]^{1/2}.
	\]
	This completes the proof.
\end{proof}

To investigate spectral redundancy, we determine the conditions under which two \( G \)-graphs have equal spectral radius. As we observe, the set \(\mathcal{S}(G_{p,q})\) encompasses three types of graphs: those of the form \(G_{p,q}\) where \(p, q \neq 0\), those where \(p = 0\) (i.e., \(G_{0,q}\)), and those where \(q = 0\) (i.e., star graphs \(G_{p,0}\)).

Now, we present a lemma that establishes key relationships between the spectral radius of a graph and its induced subgraph. This lemma serves as an important tool for the subsequent analysis.
\begin{lemma}\label{ind}\em\cite{induced_lemma}
Let \( H \) be a proper induced subgraph of \( G \). Then \(\rho(H) < \rho(G)\).
\end{lemma}

It is evident from Lemma \ref{ind} that two non-isomorphic \( G \)-graphs with zero pendent vertices cannot have equal spectral radius, as in such cases, one graph would necessarily be isomorphic to an induced subgraph of the other. Similarly, two star graphs cannot have equal spectral radius. From equation~(\ref{spectral_radius}), one can see that two non-isomorphic garlic graphs, say \( G_{p_1,q_1} \) and \( G_{p_2,q_2} \), have the same spectral radius if and only if
\[
p_1 + 2q_1 + \sqrt{p_1^2 + 4q_1^2} = p_2 + 2q_2 + \sqrt{p_2^2 + 4q_2^2}.
\]
Our goal is to determine conditions on \( p_i, q_i \) for \( i=1,2 \) under which two non-isomorphic garlic graph \( G_{p_1,q_1} \) and \( G_{p_2,q_2} \) have the same spectral radius. We begin by analyzing two non-isomorphic \( G \)-graphs with \( p_i \geq 1 \) and \( q_i \geq 1 \) for \( i=1,2 \) that share the same spectral radius. The following theorem formalizes these conditions.
\begin{theorem}\label{main_theorem}
	Let \( (p_1, q_1) \) and \( (p_2, q_2) \) be distinct pairs with \( p_i \geq 1, q_i \geq 1 \,\, \text{for}\, i=1,2\). Then, the graphs \( G_{p_1, q_1} \) and \( G_{p_2, q_2} \) have the same spectral radius if and only if one of the following conditions holds:
	\begin{itemize}
		\item [(i)] \( r_1 \coloneq \sqrt{p_1^2 + 4q_1^2} \) and \( r_2 \coloneq \sqrt{p_2^2 + 4q_2^2} \) are distinct integers such that
		\[
		p_1 + 2q_1 + r_1 = p_2 + 2q_2 + r_2.
		\]
		\item[(ii)] \( p_1 = 2q_2 \) and \( p_2 = 2q_1 \).
	\end{itemize}
	In the first case, the common spectral radius is given by
	\[
	\rho = \sqrt{\frac{p_1 + 2q_1 + r_1}{2}} = \sqrt{\frac{p_2 + 2q_2 + r_2}{2}}.
	\]
	In the second case, the common spectral radius takes the form
	\[
	\rho = \left[ q_1 + q_2 + \sqrt{q_1^2 + q_2^2} \right]^{1/2}.
	\]
\end{theorem}



\begin{proof}
Let \( G_{p_1,q_1} \) and \( G_{p_2,q_2} \) be two non-isomorphic graphs in \( \mathcal{G}\) with the same spectral radius \( \rho \). We know that the spectral radius of a garlic graph \( G_{p,q} \) is the largest root of equation \eqref{char_eq}. Therefore, for each graph \( G_{p_i, q_i} \),  \( i = 1, 2 \), we have
\begin{equation}\label{13}
\rho^4 - (p_i + 2q_i) \rho^2 + p_i q_i = 0.
\end{equation}
Since both graphs have the same spectral radius \( \rho \), therefore, we must have
\[
\rho^4 - (p_1 + 2q_1) \rho^2 + p_1 q_1 = \rho^4 - (p_2 + 2q_2) \rho^2 + p_2 q_2.
\]
By rearranging the coefficients of like terms, we obtain
\begin{equation}\label{20}
\rho^2 \left( (p_1 + 2q_1) - (p_2 + 2q_2) \right) = p_1 q_1 - p_2 q_2.
\end{equation}
Now, two cases arise.
\vspace{0.2cm}\\
\textbf{Case 1.} \( (p_1 + 2q_1) \neq (p_2 + 2q_2) \)
\vspace{0.2cm}\\
If \( (p_1 + 2q_1) \neq (p_2 + 2q_2) \), then equation (\ref{20}) can be solved for \( \rho^2 \)
\[
\rho^2 = \frac{p_1 q_1 - p_2 q_2}{(p_1 + 2q_1) - (p_2 + 2q_2)}.
\]
Since \( \rho^2 \) is a rational root of a monic polynomial $P(x)=x^2 - (p_i + 2q_i) x + p_i q_i$, therefore it must be an integer. Using equation (\ref{13}), we also have
\begin{equation}\label{sp_rad}
    \rho^2 = \frac{p_i + 2q_i + \sqrt{p_i^2 + 4q_i^2}}{2}, \,\, \text{for}\,\, i = 1, 2.
\end{equation}
For \( \rho^2 \) to be an integer, the term under the square root must be a perfect square. This implies
\[
p_1^2 + 4q_1^2 = r_1^2, \quad \text{and} \quad p_2^2 + 4q_2^2 = r_2^2
\]
for some integers \( r_1 \) and \( r_2 \). Furthermore, since both graphs share the same spectral radius, equation \eqref{sp_rad} gives the condition
\[
p_1 + 2q_1 + r_1 = p_2 + 2q_2 + r_2.
\]
This establishes that condition $(i)$ is necessary in this case.

Conversely, assume that \( p_1 + 2q_1 + r_1 = p_2 + 2q_2 + r_2 \), where $r_1$ and $r_2$ are two integers satisfying $r_i= \sqrt{p_i^2+q_i^2}$,  $i= 1,2$. Then, by equation \eqref{sp_rad}, the graphs \( G_{p_1,q_1} \) and \( G_{p_2,q_2} \) indeed have the same spectral radius.
\vspace{0.2cm}\\
\textbf{Case 2.} \( (p_1 + 2q_1) = (p_2 + 2q_2) \)
\vspace{0.2cm}\\
In this case, equation \eqref{20} reduces to
\(
p_1 q_1 = p_2 q_2.
\)
To explore this further, consider the equation \( p_1 + 2q_1 = p_2 + 2q_2 \). Squaring both sides and rearranging terms gives
\[
(p_1 + 2q_1)^2 - 8p_1 q_1 = (p_2 + 2q_2)^2 - 8p_2 q_2 \implies (p_1 - 2q_1)^2 = (p_2 - 2q_2)^2.
\]
Taking square roots yields two possibilities: $(p_1 - 2q_1) = \pm (p_2 - 2q_2)$\\
When taking the positive sign, we obtain \( p_1 = p_2 \) and \( q_1 = q_2 \), which contradicts the assumption that the graphs are non-isomorphic. When taking the negative sign, we get
\[
p_1 = 2q_2 \quad \text{and} \quad p_2 = 2q_1.
\]
Conversely, assume \( p_1 = 2q_2 \) and \( p_2 = 2q_1 \). By equation \eqref{sp_rad}, graphs \( G_{p_1,q_1}^1 \) and \( G_{p_2,q_2}^2 \) have the same spectral radius, which is given by
\[
\rho = \left[q_1 + q_2 + \sqrt{q_1^2 + q_2^2}\right]^{1/2}.
\]
Both cases establish the conditions under which two non-isomorphic graphs in \( \mathcal{G}((p,q)) \) have the same spectral radius. This completes the proof.
\end{proof}
In Theorem \ref{main_theorem} $(i)$, the triplets \( (p_1, 2q_1, r_1) \) and \( (p_2, 2q_2, r_2) \) form distinct Pythagorean triplets (PTs) with the same perimeter. Here, the perimeter of a PT refers to the perimeter of the right-angled triangle defined by the sides specified in the PT. Notably, there are infinitely many pairs of PTs with the same perimeter, as scaling any such pair by a positive integer \( c \) results in another pair with the perimeter scaled by \( c \). This is illustrated by the following example.\\

\noindent \textbf{Example.} Consider the Pythagorean triplets \((10, 24, 26)\) and \((15, 20, 25)\), both having perimeter of 60. This is the smallest perimeter shared by two distinct right-angled triangles with integer sides. The graphs \( G_{10,12} \) and \( G_{15,10} \), corresponding to these triplets, both share the same spectral radius, which is \(\sqrt{60}\). Furthermore, any non-zero positive integer multiple of these triplets, \((10c, 24c, 26c)\) and \((15c, 20c, 25c)\), will also have the same perimeter, given by \(60c\), where \(c\) is a non-zero positive integer.\\

\noindent A few more pairs of Pythagorean triplets with same perimeter, are shown in Table \ref{tab:tab1}.
\begin{table}[ht]
\centering
\begin{tabular}{|p{2cm}|p{2cm}|p{2cm}|p{2cm}|p{2cm}|p{2cm}|c|c|c|}
\hline
\centering\cellcolor{gray!25} \textbf{\bf{Triplets}} & $(21, 28, 35)$
   $(12, 35, 37)$ & $(15, 36, 39)$  $(9, 40, 41)$ & $(30, 40, 50)$  $(20, 48, 52)$ & $(33, 44, 55)$ $(11, 60, 61)$ \\ \hline
\centering \cellcolor{gray!25} \textbf{Perimeter} & \qquad 84 & \qquad 90 &
\qquad 120 & \qquad 132 \\ \hline
\end{tabular}
\caption{ Triplets with common perimeter}
\label{tab:tab1}
\end{table}

Next, we turn our attention to the conditions under which \( G_{p_1, q_1} \) and \( G_{p_2, q_2} \), where at least one of \( p_1, p_2, q_1, \) or \( q_2 \) is zero, share the same spectral radius.
The following theorems formally establish these conditions.
\begin{theorem}\label{main_thm2}
The star graph \( S_n \) and the garlic graph \( G_{p,q} \) $(q\neq 0)$ have the same spectral radius if and only if one of the following holds
\begin{enumerate}
    \item [(i)] For \( p \neq 0 \), \( n = \frac{p + 2q + r}{2} + 1 \), where, $r\coloneq \sqrt{p^2+4q^2}$ is an positive integer.
    \item [(ii)]For \( p = 0 \), \( n = 2q + 1 \).
\end{enumerate}
\end{theorem}
\begin{proof}
We consider two graphs, a star graph \( S_n \) and a garlic graph \( G_{p,q} \). Assume that they have a common spectral radius \( \rho \). First, consider the case where \( p \neq 0 \). It is well known that the spectral radius \( \rho \) of the star graph \( S_n \) is given by
\[
\rho^2 = n - 1.
\]
The spectral radius of the garlic graph \( G_{p,q} \) is given by
\[
\rho^2 = \frac{p + 2q + \sqrt{p^2 + 4q^2}}{2}.
\]
Equating the two expressions for \( \rho^2 \), we obtain
\begin{equation}\label{eq1}
n - 1 = \frac{p + 2q + \sqrt{p^2 + 4q^2}}{2}.
\end{equation}
Since \( n - 1 \) is an integer, the right-hand side of equation \eqref{eq1} must also be an integer. This requires \( \sqrt{p^2 + 4q^2} \) to be an integer. Let \( \sqrt{p^2 + 4q^2} = r \), where \( r \) is an integer. Substituting \( r = \sqrt{p^2 + 4q^2} \) into equation \eqref{eq1}, we can solve for \( n \) to obtain
\[ n = \frac{p + 2q + r}{2} + 1. \]
Conversely, if \( n = \frac{p + 2q + r}{2} + 1 \), where $p$, $q$ and $r=\sqrt{p^2+4q^2}$ are positive integers then, the spectral radius of \( G_{p,q} \) is given by
\[
\rho = \sqrt{\frac{p + 2q + r}{2}}.
\]
For \( n = \frac{p + 2q + r}{2} + 1 \), the spectral radius of \( S_n \) is
\[
\rho = \sqrt{n - 1} = \sqrt{\frac{p + 2q + r}{2}},
\]
matching that of \( G_{p,q} \). Thus, the two graphs have the same spectral radius.\\

Now, we consider the case when \( p = 0 \). In this scenario, equation \eqref{eq1} reduces to
\[
n - 1 = 2q, \quad \text{so} \quad n = 2q + 1.
\]
Conversely, assume that \( n = 2q + 1 \). The spectral radius of \( S_n \) is
\[
\rho = \sqrt{n - 1} = \sqrt{2q}.
\]
For the garlic graph \( G_{p,q} \) with \( p = 0 \), its spectral radius is
\[
\rho = \sqrt{\frac{p + 2q + \sqrt{p^2 + 4q^2}}{2}} = \sqrt{\frac{0 + 2q + \sqrt{0 + 4q^2}}{2}} = \sqrt{2q}.
\]
Thus, the spectral radius of \( S_n \) matches that of \( G_{0,q} \) when \( n = 2q + 1 \). This completes the proof. \end{proof}

\noindent \textbf{Remark.} The condition \( q \neq 0 \) is necessary because if \( q = 0 \), one graph would become a proper induced subgraph of the other, which contradicts the possibility of equal spectral radius.
\begin{theorem}\label{main_thm3}
	Let \( p_2 \geq 1 \), \( q_1 \geq 1 \), and \( q_2 \geq 1 \). Then the graphs \( G_{0,q_1} \) and \( G_{p_2,q_2} \) have the same spectral radius if and only if  \(	r_2 \coloneqq \sqrt{p_2^2 + 4q_2^2}\) is an integer such that
	\(	4q_1 = p_2 + 2q_2 + r_2.\)
\end{theorem}
\begin{proof}
The proof follows a similar structure to that of the preceding theorem and is therefore omitted.
\end{proof}
\begin{theorem}\label{thm27}
Let \( (p, q, r) \) or \( (q, p, r) \) be a Pythagorean triplet , then there exist precisely four graphs, \( S_{p + q + r + 1} \) (or \( G_{p+q+r, 0} \)), \( G_{0, \frac{p+q+r}{2}} \), \( G_{2q, p} \), and \( G_{2p, q} \), all with spectral radius \( \sqrt{p + q + r} \).
\end{theorem}
\begin{proof}
Let \( (p, q, r) \) be a Pythagorean triplet. Consider the graphs \( G^1_{p_1, q_1} \) and \( G^2_{p_2, q_2} \), where \( p_1 = 2q \), \( q_1 = p \), \( p_2 = 2p \), and \( q_2 = q \). These parameter choices satisfy the conditions of Theorem \ref{main_theorem} $(ii)$, ensuring that both \( G^1_{p_1, q_1} \) and \( G^2_{p_2, q_2} \) have the same spectral radius, given by \( \sqrt{p + q + r} \).

Furthermore, by Theorem \ref{main_thm2} and Theorem \ref{main_thm3} , this spectral radius is also equal to that of the star graph \( S_{p + q + r + 1} \) and the garlic graph \( G_{0, (p + q + r)/2} \). Hence, the graphs \( G_{2q, p} \), \( G_{2p, q} \), \( G_{0, (p + q + r)/2} \), and \( S_{p + q + r + 1} \) all share the spectral radius \( \sqrt{p + q + r} \).

To demonstrate that no other garlic graph, non-isomorphic to the four given graphs, corresponds to the Pythagorean triplet \( (p, q, r) \), it is necessary to establish that any garlic graph with a spectral radius \( \sqrt{p + q + r} \), other than the four previously identified, must correspond to a distinct Pythagorean triplet. The spectral radius of the graphs \( S_{p + q + r + 1} \) and \( G_{0, \frac{p + q + r}{2}} \) depend on the perimeter of the Pythagorean triplet \((p, q, r)\), meaning no additional non-isomorphic graphs can be formed in this category with the spectral radius \( \sqrt{p + q + r} \).

Now, assume there exists a graph \( G_{a, b} \) with \( a, b \neq 0 \) that is non-isomorphic to the four stated graphs but shares the same spectral radius. According to Theorem \ref{main_thm3}, for \( G_{a, b} \) to have the same spectral radius as \( G_{0, \frac{p + q + r}{2}} \), there must be a Pythagorean triplet \((a, 2b, c)\), where \( c = \sqrt{a^2 + (2b)^2} \). Additionally, the condition \(\frac{p + q + r}{2} = \frac{a + 2b + c}{4}\) must hold. Simplifying this equation leads to
\[
p + q + r = \frac{a + 2b + c}{2},
\]
indicating that the perimeter of one Pythagorean triplet is half that of the other. This implies the existence of another Pythagorean triplet \((a, 2b, c)\).

Furthermore, it is easy to verify that the triplet \( (q, p, r) \) yields the same set of four graphs, thereby completing the proof.
\end{proof}

Each Pythagorean triplet \( (p, q, r) \) corresponds to four graphs that share the same spectral radius, \( \sqrt{p+q+r} \). However, as discussed earlier, multiple Pythagorean triplets can have the same perimeter. This means that more than four graphs can share the same spectral radius \( \sqrt{p+q+r} \). Consider the following example for a more detailed illustration.\\

\noindent \textbf{Example.} The Pythagorean triplets \((15, 20, 25)\) and \((10, 24, 26)\) both have a sum of 60. As a result, the graphs \(G_{40, 15}\), \(G_{30, 20}\), \(G_{48, 10}\), \(G_{20, 24}\), $G_{0,30}$ and \(S_{61}\) all share the same spectral radius \(\sqrt{60}\).

For the triplet \((15, 20, 25)\), the perimeter is \(15 + 20 + 25 = 60\), and the corresponding graphs are \(G_{40, 15}\) and \(G_{30, 20}\). Similarly, for the triplet \((10, 24, 26)\), the perimeter is \(10 + 24 + 26 = 60\), and the corresponding graphs are \(G_{48, 10}\) and \(G_{20, 24}\).

Both sets of graphs share the same spectral radius because they have the same perimeter. Additionally, the star graph \(S_{61}\) (where \(n = p + q + r + 1\)) and $G_{0,30}$ also shares the same spectral radius, \(\sqrt{60}\).\\

Table \ref{tab:tab2} displays multiple Pythagorean triplets that share a common perimeter. This raises an important question: What is the maximum number of Pythagorean triplets that can have the same perimeter, and consequently, how many graphs can share the same spectral radius? To address this, we observe that for any positive integer \( m \), there exists an integer \( P \) such that \( m \) distinct Pythagorean triplets share the same perimeter \( P \) and at least $2m+1$ non-isomorphic garlic graphs have common spectral radius. The following theorem proves that it is indeed possible for \( m \) Pythagorean triplets to have the same perimeter and \( 2m + 1 \) garlic graphs to share the same spectral radius.
\begin{table}[ht]
\centering
\begin{tabular}{|c|p{8cm}|c|}
\hline
\cellcolor{gray!25} \textbf{No. of Triplets (m)} & \cellcolor{gray!25} \textbf{ \qquad \qquad \qquad \qquad Triplets} & \cellcolor{gray!25} \textbf{Perimeter (P)} \\ \hline
2 & $(15, 20, 25)$, $(10, 24, 26)$ & 60 \\ \hline
3 & $(30, 40, 50)$, $ (20, 48, 52)$, $(24, 45, 51)$ & 120 \\ \hline
4 & $(60, 80, 100)$,\,$(40, 96, 104)$,\,$(48, 90, 102)$,
$(15, 112, 113)$ & 240 \\ \hline
5 & $(105, 140, 175)$, $(70, 168, 182)$, $(120, 126, 174)$, $(60, 175, 185)$, $(28, 195, 197)$ & 420 \\ \hline
\end{tabular}
\caption{Triplets with common perimeter.}
\label{tab:tab2}
\end{table}
\begin{theorem} For any positive integer $m$ there exists at least $2m+1$ graphs in $\mathcal{G}$ having common spectral radius.   
\end{theorem}
\begin{proof}
Let \( m \) be a positive integer, and consider the primitive Pythagorean triplets \( T_1, T_2, \ldots, T_m \) obtained from the general formula \( (a^2 - b^2, 2ab, a^2 + b^2) \), where \( a \in \{ p_1, p_2, \ldots, p_m \} \) and \( b = 2 \), with \( p_i \) representing the distinct odd prime numbers. For each \( T_i \), we denote the perimeter by \( P_i \). Let \( P \) be the least common multiple (l.c.m.) of the perimeters \( P_1, P_2, \ldots, P_m \). We can then construct scaled Pythagorean triplets \( \frac{P T_1}{P_1}, \frac{P T_2}{P_2}, \ldots, \frac{P T_m}{P_m} \) such that each triplet shares the same perimeter $P$.

To show that these scaled triplets are distinct, assume that there exist two triplets \( \frac{P T_j}{P_j} \) and \( \frac{P T_k}{P_k} \) that are equal for \( j \neq k \). This implies
\begin{align*}
\frac{P}{2p_j(p_j + 2)}(p_j^2 - 4, 4p_j, p_j^2 + 4) &= \frac{P}{2p_k(p_k + 2)}(p_k^2 - 4, 4p_k, p_k^2 + 4), \\
\Rightarrow \quad 2p_k(p_k + 2)(p_j^2 - 4, 4p_j, p_j^2 + 4) &= 2p_j(p_j + 2)(p_k^2 - 4, 4p_k, p_k^2 + 4).
\end{align*}
For these two scaled triplets to be equal, the third components must also be equal
\[
2p_k(p_k + 2)(p_j^2 + 4) = 2p_j(p_j + 2)(p_k^2 + 4).
\]
Simplifying this equation gives
\[
p_k | (p_j + 2)(p_k^2 + 4) \implies p_k | (p_j + 2),
\]
and in similar way we obtain
\[
p_j | (p_k + 2).
\]
Assume without loss of generality that \( p_k > p_j \). Since \( p_k \mid (p_j + 2) \), we must have \( p_j + 2 \geq p_k \). However, this leads to the inequality \( p_j + 2 \geq p_k > p_j \), since the difference between two odd primes must be at least 2. This can only be true if
\[
p_k = p_j + 2.
\]
Next, from \( p_j \mid (p_k + 2) \), we have
\[
p_j \mid (p_j + 4) \Rightarrow  p_j \mid 4.
\]
However, since \( p_j \) is an odd prime, this leads to a contradiction. Thus, no distinct primes \( p_j \) and \( p_k \) can satisfy the original assumption that \( \frac{P T_j}{P_j} \) and \( \frac{P T_k}{P_k} \) are equal. Therefore, all scaled Pythagorean triplets generated by the described process are distinct. According to Theorem \ref{thm27}, each triplet corresponds to exactly four graphs with the same spectral radius. Among these, two graphs are isomorphic to \( S_{P + 1} \) and \( G_{0, P/2} \) in every quadruple due to the equal perimeter of all triplets. For the remaining two graphs in each quadruple, we consider the set of \( 2m \) graphs formed by their collection and demonstrate that each graph in this set is distinct and non-isomorphic to all others.

Let \( G_1 \) and \( G_2 \) be two graphs in this set. They either correspond to the same Pythagorean triplet or to distinct Pythagorean triplets. If \( G_1 \) and \( G_2 \) correspond to the same Pythagorean triplet \((a, b, c)\), then \( G_1 = G_{2b, a} \) and \( G_2 = G_{2a, b} \). For \( G_1 \) and \( G_2 \) to be isomorphic, it must hold that \( 2b = 2a \), which implies \( a = b \). However, this contradicts the fact that no two sides of a Pythagorean triplet can be equal. Therefore, \( G_1 \) and \( G_2 \) cannot be isomorphic.

If \( G_1 \) and \( G_2 \) correspond to distinct triplets \((a, b, c)\) and \((x, y, z)\), then \( G_1 \in \{ G_{2b, a}, G_{2a, b} \} \) and \( G_2 \in \{ G_{2y, x}, G_{2x, y} \} \), resulting in four possibilities: \( G_{2b, a} \cong G_{2x, y} \), \( G_{2b, a} \cong G_{2y, x} \), \( G_{2a, b} \cong G_{2x, y} \), or \( G_{2a, b} \cong G_{2y, x} \). Assuming any of these cases, e.g., \( G_{2b, a} \cong G_{2x, y} \), implies \( 2b = 2x \) and \( a = y \), leading to \((a, b, c) = (y, x, z)\). Consequently, \( G_1 \) and \( G_2 \) correspond to \((y, x, z)\) and \((x, y, z)\), respectively.

Next, we show that no such pair of triplets exists among the generated PTs by proving that if \((x, y, z)\) is among the generated Pythagorean triplets (PTs), then \((y, x, z)\) cannot be. To do this, suppose both \((x, y, z)\) and \((y, x, z)\) are among the generated PTs. Then there exist two primes \( p_i \) and \( p_j \) such that
\[
(x, y, z) = \frac{P}{2p_j(p_j + 2)}(p_j^2 - 4, 4p_j, p_j^2 + 4)
\quad \text{and} \quad
(y, x, z) = \frac{P}{2p_k(p_k + 2)}(p_k^2 - 4, 4p_k, p_k^2 + 4).
\]
Equating the values of $x$ from both PTs implies
\[
\frac{P (p_i^2 - 4)}{2p_i(p_i + 2)} = \frac{4Pp_j}{2p_j(p_j + 2)}.
\]
Simplifying this expression yields
\[
(p_i^2 - 4)(p_j + 2) = 4p_i(p_i + 2).
\]

Since \( p_i \) and \( p_j \) are odd primes, the left-hand side of this equation is odd, while the right-hand side is even, leading to a contradiction. This contradiction shows that such a pair of distinct triplets cannot exist. Similar contradictions arise for the other cases, ensuring that \( G_1 \) and \( G_2 \) are not isomorphic.

Thus, among the \( 4m \) graphs, at least \( 2m + 2 \) graphs, including the star graph, are pairwise non-isomorphic. However, since the star graph does not belong to the family \(\mathcal{G}\), it is excluded. Therefore, there are at least \( 2m + 1 \) non-isomorphic graphs with a common spectral radius \( \sqrt{P} \).
\end{proof}

\section{Spectral redundancy index}

Determining an explicit formula of \( c(G_{p, q}) \) for the graph $G_{p, q}$ in \( \mathcal{G} \) is not feasible due to the difficulty of explicitly counting Pythagorean triplets. However, upper and lower bounds for \( c(G_{p,q}) \) can be established. To achieve this, we first establish several intermediate results that will be required. For this purpose, we define two sets of unordered pairs of graphs as follows.

\[
S_1 = \left\{ \left( G_{p_1, q_1}, G_{p_2, q_2} \right) \,|\, G_{p_1, q_1}, G_{p_2, q_2} \in \mathcal{S}(G_{p, q}), \, p_1 = 2q_2, \, p_2 = 2q_1 \right\},
\]

\[
S_2 = \left\{ \left( S_{n}, G_{0, q_1} \right) \,|\, S_{n}, G_{0, q_1} \in \mathcal{S}(G_{p, q}), \, n = 2q_1 + 1 \right\}.
\]

By Theorems \ref{main_theorem} and \ref{main_thm2}, both graphs in each pair within the sets \( S_1 \) and \( S_2 \) must have equal spectral radius. The following lemma investigates whether two pairs from these sets share a common spectral radius.

\begin{lemma}\label{lem:s1s2}
Let \( S_1 \) and \( S_2 \) be the sets defined above. The following properties hold.

\begin{enumerate}
    \item[(a)] Let \( (G_{p_1, q_1}, G_{p_2, q_2}) \) and \( (G_{p'_1, q'_1}, G_{p'_2, q'_2}) \)  \( \in  S_1 \). These pairs share a common spectral radius if and only if \( (q_1, q_2, r) \) and \( (q'_1, q'_2, r') \) are Pythagorean triplets for some integer $r, r'$, such that \( q_1 + q_2 + r = q'_1 + q'_2 + r' \).

    \item[(b)] No two pairs of graphs in the set \( S_2 \) share a common spectral radius.

    \item[(c)] Let \( (G_{p_1, q_1}, G_{p_2, q_2}) \in S_1 \) and \( (S_{n}, G_{0, q_3}) \in S_2 \). These two pairs share the same spectral radius if and only if \( (q_1, q_2, r) \) forms a Pythagorean triplet, where \( q_1^2 + q_2^2 = r^2 \) for some integer \( r \), and \( q_3 \) satisfies \( 2q_3 = q_1 + q_2 + r \).
\end{enumerate}
\end{lemma}

\begin{proof}
Suppose \( (G_{p_1, q_1}, G_{p_2, q_2}) \) and \( (G_{p'_1, q'_1}, G_{p'_2, q'_2}) \) are two pairs of graphs in \( \mathcal{S}_2 \) that share a common spectral radius. By Theorem \ref{main_theorem}, we have
\[
q_1 + q_2 + \sqrt{q_1^2 + q_2^2} = q'_1 + q'_2 + \sqrt{q_1'^2 + q_2'^2}.
\]
Let \( A_1 = q_1 + q_2 \), \( A_2 = q'_1 + q'_2 \), \( B_1 = q_1^2 + q_2^2 \), and \( B_2 = q_1'^2 + q_2'^2 \). This simplifies to
\[
A_1 + \sqrt{B_1} = A_2 + \sqrt{B_2}.
\]
If \( B_1 \) and \( B_2 \) are not perfect squares, \( \sqrt{B_1} \) and \( \sqrt{B_2} \) are irrational. For equality, both \( A_1 = A_2 \) and \( B_1 = B_2 \) must hold. Thus
\[
q_1 + q_2 = q'_1 + q'_2 \quad \text{and} \quad q_1^2 + q_2^2 = q_1'^2 + q_2'^2.
\]
Expanding \( (q_1 + q_2)^2 = (q'_1 + q'_2)^2 \) results in \( q_1q_2 = q'_1q'_2 \). Combining this with \( q_1^2 + q_2^2 = q_1'^2 + q_2'^2 \), we find \( (q_2 - q_1)^2 = (q'_2 - q'_1)^2 \). Without loss of generality, assuming \( q_2 \geq q_1 \) and \( q'_2 \geq q'_1 \), it follows that \( q_2 - q_1 = q'_2 - q'_1 \). Solving \( q_1 + q_2 = q'_1 + q'_2 \) and \( q_2 - q_1 = q'_2 - q'_1 \) gives \( q_1 = q'_1 \) and \( q_2 = q'_2 \), proving that the pairs are identical.

Therefore, for distinct pairs to share the same spectral radius, \( B_1 \) and \( B_2 \) must be perfect squares. Consequently, there exist integers \( r \) and \( r' \) such that \( r = \sqrt{q_1^2 + q_2^2} \) and \( r' = \sqrt{q_1'^2 + q_2'^2} \), forming Pythagorean triplets \( (q_1, q_2, r) \) and \( (q'_1, q'_2, r') \) with
\[
q_1 + q_2 + r = q'_1 + q'_2 + r'.
\]
Conversely, if \( (q_1, q_2, r) \) and \( (q'_1, q'_2, r') \) are two Pythagorean triplets with
\[
q_1 + q_2 + \sqrt{q_1^2 + q_2^2} = q'_1 + q'_2 + \sqrt{q_1'^2 + q_2'^2}
\]
implies the graphs \( (G_{p_1, q_1}, G_{p_2, q_2}) \) and \( (G_{p'_1, q'_1}, G_{p'_2, q'_2}) \) have the same spectral radius, completing the proof for $(a)$.

Let \( (S_{n_1}, G_{0, q_1}), (S_{n_2}, G_{0, q_2}) \in S_2 \). By Theorem \ref{main_thm2}, their spectral radius are \( \sqrt{2q_1} \) and \( \sqrt{2q_2} \), respectively. If these pairs share a common spectral radius, then \( \sqrt{2q_1} = \sqrt{2q_2} \), which implies \( q_1 = q_2 \). Consequently, since \( n_1 = 2q_1 + 1 \) and \( n_2 = 2q_2 + 1 \), it follows that \( n_1 = n_2 \). Therefore, the pairs are isomorphic. This completes the proof for $(b)$.

Let \( (G_{p_1, q_1}, G_{p_2, q_2}) \in S_1 \) and \( (S_{n}, G_{0, q_3}) \in S_2 \) share the same spectral radius. By Theorems \ref{main_theorem} and \ref{main_thm2}, this implies
\[
2q_3 = q_1 + q_2 + \sqrt{q_1^2 + q_2^2}.
\]
The left-hand side is an integer, so the right-hand side must also be an integer. This requires \( \sqrt{q_1^2 + q_2^2} \) to be an integer, say \( r \), such that \( q_1^2 + q_2^2 = r^2 \). Hence, \( (q_1, q_2, r) \) forms a Pythagorean triplet with condition $2q_3 = q_1 + q_2 + r.$

Conversely, suppose \( (q_1, q_2, r) \) forms a Pythagorean triplet, i.e., \( q_1^2 + q_2^2 = r^2 \), where \( r \) is an integer. Define \( q_3 \) such that
\[
2q_3 = q_1 + q_2 + r.
\]
 The spectral radius of \( (S_{n}, G_{0, q_3}) \), given by \( \sqrt{2q_3} \), matches the spectral radius of \( (G_{p_1, q_1}, G_{p_2, q_2}) \), given by \( \sqrt{q_1 + q_2 + r} \). This completes the proof of the lemma.
\end{proof}

The following theorem provides a lower bound for \( c(G_{p,q}) \).
\begin{theorem}\label{lower_bound}
    Let \( G_{p,q} \) be a garlic graph. Then \(  b(G_{p,q}) - \frac{(k-1)(k+2))}{2} \leq c(G_{p,q}) \), where \( k = \min \left\{ q, \left\lfloor \frac{p}{2} \right\rfloor \right\} \). Equality holds if and only if \((p, q) \in \{1, 2\} \times \mathbb{N}\) or \((p, q) \in \mathbb{N} \times \{1\}\)
\end{theorem}
\begin{proof}Let \( G_{p,q} \) be a garlic graph, and let \( \mathcal{S}(G_{p,q}) \) denote the set of all induced non-isomorphic graphs of \( G_{p,q} \). It has been observed that the presence of Pythagorean triplets plays a significant role in contributing to the redundancy of graphs with equal spectral radius in \( \mathcal{S}(G_{p,q}) \). To establish a lower bound for \( c(G_{p,q}) \), we consider the scenario where no Pythagorean triplet exists.

In such cases, graphs with equal spectral radius must satisfy the conditions outlined in Theorem \ref{main_theorem} (ii) or Theorem \ref{main_thm2} (ii). Furthermore, all pairs of graphs with equal spectral radius must either belong to \( S_1 \) or \( S_2 \). Lemma \ref{lem:s1s2} confirms that, in the absence of Pythagorean triplets, no pairs in \( S_1 \) and \( S_2 \) share the same spectral radius. Next, we proceed to count the pairs in \( S_1 \) and \( S_2 \) under these conditions.

We consider two cases based on the relationship between \( p \) and \( q \). First, assume that \( p \geq 2q \).
For each pair \( (q_1, q_2) \), where \( 1 \leq q_1, q_2 \leq q \), there exist corresponding values \( p_1 = 2q_1 \) and \( p_2 = 2q_2 \) such that the graphs \( G_{p_1, q_1} \) and \( G_{p_2, q_2} \) are induced subgraphs of \( G_{p,q} \). By part $(ii)$ of Theorem \ref{main_theorem}, these graphs share the same spectral radius. The total number of such pairs \( (q_1, q_2) \) is \( \frac{q(q-1)}{2} \). Additionally, for each \( 1 \leq q_1 \leq q \), Theorem \ref{main_thm2} $(ii)$ establishes that the graphs \( G_{0, q_1} \) and \( G_{2q_1, 0} \) also share the same spectral radius. This contributes \( q \) additional pairs. However, the graphs \( G_{2, 0} \) and \( G_{0, 1} \) are isomorphic to each other and therefore are excluded.

Thus, the total number of pairs of induced subgraphs sharing a common spectral radius is
\[
\frac{q(q-1)}{2} + (q-1) = \frac{(q-1)(q+2)}{2}.
\]

Similarly, consider the case when \( 2q > p \). For every pair of even integers \( (p_1, p_2) \), where \( 1 < p_1, p_2 \leq p \), there exist corresponding values \( q_1 = \frac{p_1}{2} \) and \( q_2 = \frac{p_2}{2} \) such that \( G_{p_1, q_1} \) and \( G_{p_2, q_2} \) are induced subgraphs of \( G_{p,q} \). By Theorem \ref{main_theorem}, these graphs share the same spectral radius. The total number of such pairs \( (p_1, p_2) \) is
\(\lfloor \frac{p}{2} \rfloor \left( \lfloor \frac{p}{2} \rfloor - 1 \right)/2.\) Furthermore, for each even integer \( 1 \leq p_1 \leq p \), Theorem \ref{main_thm2} $(ii)$ guarantees that the graphs \( G_{0, \frac{p_1}{2}} \) and \( G_{p_1, 0} \) share the same spectral radius, contributing \( \lfloor \frac{p}{2} \rfloor \) additional pairs. As before, the graphs \( G_{2, 0} \) and \( G_{0, 1} \) are isomorphic to each other and are therefore excluded.

Hence, the total number of pairs of induced subgraphs sharing a common spectral radius is
\[
\frac{\lfloor \frac{p}{2} \rfloor \left( \lfloor \frac{p}{2} \rfloor - 1 \right)}{2} + \left\lfloor \frac{p}{2} \right\rfloor -1= \frac{\left(\lfloor \frac{p}{2}\rfloor-1  \right)  \left( \lfloor \frac{p}{2} \rfloor + 2 \right)}{2}.
\]
Combining the results from the two cases, we conclude that there are \( \frac{(k-1)(k+2)}{2} \) pairs with a common spectral radius in $\mathcal{S}(G_{p,q})$, where \( k = \min \left\{ q, \left\lfloor \frac{p}{2} \right\rfloor \right\} \). Consequently, the following expression for \( c(G_{p,q}) \) is obtained

\[
c(G_{p,q}) = b(G_{p,q}) - \frac{(k-1)(k+2)}{2},
\]

This equation was derived under the assumption that no Pythagorean triplets exist, and therefore, we did not account for any repetition in spectral radius caused by Pythagorean triplets. However, in reality, Pythagorean triplets do exist. Thus, we conclude that in all cases, the inequality holds as
\[ b(G_{p,q}) - \frac{k(k+1)}{2} \leq c(G_{p,q}),\]

Since, there are no triplets with side lengths less than or equal to 2, this scenario falls within the earlier assumption, so equality holds in situations where \( (p, q) \in \{1, 2\} \times \mathbb{N} \) or \( (p, q) \in \mathbb{N} \times \{1\} \).
\end{proof}
\begin{corollary}
    The infinite families \( \mathcal{G}_{1,q} \), \( \mathcal{G}_{2,q} \) and \( \mathcal{G}_{p,1} \) defined as
    \[
    \mathcal{G}_{1,q} = \{ G_{p,q} : p = 1, q \geq 1 \},  \mathcal{G}_{2,q} = \{ G_{p,q} : p = 2, q \geq 1 \},  \mathcal{G}_{p,1} = \{ G_{p,1} : p \geq 1,  q = 1 \},
    \]
    are spectrally non-redundant, that is, their spectral redundancy index is equal to \( 1 \).
\end{corollary}

\begin{proof}

For each family, we have \( k = \min\{ q, \lfloor \frac{p}{2} \rfloor \} = 1 \). According to the theorem, \( b(G) \leq c(G) \). Additionally, \( c(G) \leq b(G) \) holds. Combining these inequalities, we find that \( c(G) = b(G) \) for all \( G \) in the families \( \mathcal{G}_{1,q} \), \( \mathcal{G}_{2,q} \), or \( \mathcal{G}_{p,1} \). Notably, the result that \( \mathcal{G}_{p,1} \) is spectrally non-redundant aligns with the findings in \cite{See}.
\end{proof}

Table \ref{tab:tab3} presents the formula for \( c(G) \) for graphs \( G \) within certain infinite subfamilies of \( \mathcal{G} \).
\begin{table}[ht]
\centering
\begin{tabular}{|c|c|c|c|c|}
\hline
\cellcolor{gray!25} \textbf{Sub-family $G$} &
\( \mathcal{G}_{1,q}, \mathcal{G}_{2,q}, \mathcal{G}_{p,1} \) &
 \( \mathcal{G}_{3,q} \) &
 \( \mathcal{G}_{4,q}, \mathcal{G}_{5,q}, \mathcal{G}_{p,2} \) &
\( \mathcal{G}_{6,q}, \mathcal{G}_{p,3} \) \\ \hline
\cellcolor{gray!25} \textbf{$c(G)$} &
$b(G)$ &
$b(G)-1$ &
$b(G)-3$ &
$b(G)-7$ \\ \hline
\end{tabular}
\caption{Few sub-families $G$ with corresponding $c(G)$ }
\label{tab:tab3}
\end{table}

Spectral redundancy is closely tied to the presence of Pythagorean triplets. However, determining the exact number of Pythagorean triplets is challenging. To compute the spectral redundancy index, an approximate count of Pythagorean triplets is necessary. This approximation is provided by the following theorem. Subsequently, we examine the spectral redundancy of certain infinite subfamilies within the family of \(G\)-graphs. The following result, presented in \cite{PT}, offers an estimate for the number of Pythagorean triplets with both legs less than or equal to a given positive integer \(n\).
\begin{theorem}\em\cite{PT}\label{PT_est}
The number \(\tilde{T}(n)\) of Pythagorean triplets \((a, b, c)\) such that \(a \leq n\) and \(b \leq n\) (considering the triplet \((a, b, c)\) the same as \((b, a, c)\)) is given by
\[
\tilde{T}(n) = \frac{2 \log(1+\sqrt{2})}{\log 2} n \log n + O(n).
\]
\end{theorem}

\begin{theorem}
	For the family \( \mathcal{G} \) of all \( G \)-graphs \( G_{p,q} \), the spectral redundancy index \( r(\mathcal{G}) \) is \( \frac{4}{3}. \)
\end{theorem}
\vspace{0.01cm}
\begin{proof}
Let \( \mathcal{G} \) be the family of all \( G \)-graphs \( G_{p,q} \). For any \( G_{p,q} \in \mathcal{G} \), by Theorem \ref{lower_bound},
\[
b(G_{p,q}) - \frac{(k-1)(k+2)}{2} \leq c(G_{p,q}),
\]
where \( k = \min \left\{ q, \lfloor p/2 \rfloor \right\} \). Rearranging,
\[
\frac{b(G_{p,q})}{c(G_{p,q})} \leq 1 + \frac{(k-1)(k+2)}{2c(G_{p,q})} \leq 1 + \frac{(k-1)(k+2)}{2\left(qp + 2q + p - \frac{(k-1)(k+2)}{2}\right)}.
\]
The last term follows from \( c(G_{p,q}) \geq b(G_{p,q}) - \frac{(k-1)(k+2)}{2} \) and Lemma \ref{subgraph}. We now consider two cases:

\vspace{0.2cm}
\noindent
\textbf{Case 1:} If \( q < \lfloor p/2 \rfloor \), then \( k = q \), $q<2p$ and using \( p/2 - 1 \leq \lfloor p/2 \rfloor \leq p/2 \), we obtain
\[
\begin{aligned}
	\frac{b(G_{p,q})}{c(G_{p,q})} & \leq 1 + \frac{(q-1)(q+2)}{2\left(qp + 2q + p - \frac{(q-1)(q+2)}{2}\right)}\\
	& \leq 1 + \frac{(q-1)(q+2)}{4q^2 + 8q - (q-1)(q+2)} \leq \frac{4}{3}, \quad \text{for all } q \geq 1.
\end{aligned}
\]

\vspace{0.2cm}
\noindent
\textbf{Case 2:} If \( q \geq \lfloor p/2 \rfloor \), then \( k = \lfloor p/2 \rfloor \), $2q\geq p-2$ and using \( p/2 - 1 \leq \lfloor p/2 \rfloor \leq p/2 \), we have
\[
\begin{aligned}
	\frac{b(G_{p,q})}{c(G_{p,q})} & \leq 1 + \frac{(\lfloor p/2 \rfloor -1)(\lfloor p/2 \rfloor+2)}{2qp + 4q + 2p - (\lfloor p/2 \rfloor -1)(\lfloor p/2 \rfloor +2)} \coloneq g(p,q)\\
	& \leq 1 + \frac{( p/2 -1)( p/2+2)}{p(p-2) + 2(p-2) + 2p - ( p/2-2)( p/2 +1)}\leq \frac{4}{3}, \quad \text{for all } p \geq 2.
\end{aligned}
\]
Since, it is easy to verify that \( g(1,q) \leq 4/3 \), therefore we can conclude that
\begin{equation}\label{r_up}
	\mathfrak{r}(\mathcal{G}) = \sup_{G_{p,q} \in \mathcal{G}} \frac{b(G_{p,q})}{c(G_{p,q})} \leq \frac{4}{3}.
\end{equation}
To determine the lower bound, we consider the subfamily,
\[
\mathcal{G}_{\text{bal}} = \{ G_{2q,q} \mid q \geq 1 \},
\]
of balanced \( G \)-graphs. We refer to these graphs as "balanced" because the condition \( p = 2q \) ensures that the sum of degrees of all pendant vertices equals the sum of degrees of all vertices of degree \( 2 \).
	
	The number of pairs of graphs with the same spectral radius in \( \mathcal{S}(G_{2q,q}) \) is given by \( \frac{(q-1)(q+2)}{2} \), assuming the absence of Pythagorean triplets. This result can be derived in a manner similar to the approach used in the Theorem \ref{lower_bound}.
	
	However, Pythagorean triplets do exist and contribute to additional redundancy in \( \mathcal{S}(G_{2q,q}) \). Using Theorem \ref{PT_est}, the number of Pythagorean triplets \( \tilde{T}(q) \), where each triplet \( (a, b, c) \) satisfies \( a, b \leq q \), is approximately
	\[
	\tilde{T}(q) = \frac{2 \log(1+\sqrt{2})}{\log 2} q \log q + O(q).
	\]
	Since each Pythagorean triplet corresponds to at most four graphs with the same spectral radius, therefore, a rough upper bound for $c(G_{p,q})$ can be found as
	
	\[
	b(G_{2q, q}) - \frac{(q-1)(q+2)}{2} \leq c(G_{2q, q}) \leq b(G_{2q, q}) - \left[\frac{(q-1)(q+2)}{2} - 4\tilde{T}(q)\right].
	\]
	From these bounds, the spectral redundancy ratio satisfies
	\[
	\frac{1}{1 - \frac{(q-1)(q+2)}{2b(G_{2q, q})} + \frac{4\tilde{T}(q)}{b(G_{2q, q})}} \leq \frac{b(G_{2q, q})}{c(G_{2q, q})} \leq \frac{1}{1 - \frac{(q-1)(q+2)}{2b(G_{2q, q})}},
	\]
	where \( b(G_{2q, q}) = 2q^2 + 4q \).
	
	As \( q \to \infty \), the term \( \frac{\tilde{T}(q)}{b(G_{2q, q})} \to 0 \), implying that both sides of the inequality approach \( \frac{4}{3} \). Consequently, we obtain
	\[
	\lim_{q \to \infty} \frac{b(G_{2q, q})}{c(G_{2q, q})} = \frac{4}{3}.
	\]
	Since, we also have
	\[
	\frac{b(G_{2q, q})}{c(G_{2q, q})} \leq \frac{1}{1 - \frac{(q-1)(q+2)}{4(q^2+2q)}}=\frac{4(q^2+2q)}{3q^2+7q+2} \leq \frac{4}{3}, \quad \text{for all } q \geq 0.
	\]
	This inequality demonstrates that \( \frac{4}{3} \) is an upper bound, and since it is also attained asymptotically as \( q \to \infty \), it represents the supremum. Thus, summarizing the above, we have
    \begin{equation}\label{r_low}
    	\mathfrak{r}(\mathcal{G}_{bal}) = \underset{G \in \mathcal{G}_{bal}}{\sup} \frac{b(G)}{c(G)} = \frac{4}{3}.
    \end{equation}
    Combining equations \eqref{r_up} and \eqref{r_low}, we obtain

    \begin{equation}
    \frac{4}{3}\leq \mathfrak{r}(\mathcal{G}_{bal})\leq 	\mathfrak{r}(\mathcal{G}) = \sup_{G_{p,q} \in \mathcal{G}} \frac{b(G_{p,q})}{c(G_{p,q})} \leq \frac{4}{3}.
    \end{equation}
This completes the proof.
\end{proof}

\section{Discussion}

In this work, we investigated a family of graphs that we refer to as garlic graphs. These graphs exhibit significant spectral redundancy, a phenomenon where multiple non-isomorphic graphs share the same spectral radius. We explored the conditions under which two graphs in this family possess identical spectral radius and demonstrated a connection between this redundancy and the presence of Pythagorean triplets. Specifically, we identified how the parameters of garlic graphs relate to the sides of Pythagorean triplets, providing a combinatorial perspective on spectral properties.

An important finding about the family of \( G \)-graphs is that we can construct an arbitrarily large set of non-isomorphic graphs, all sharing the same spectral radius. Furthermore, we determine the spectral redundancy index of the family of \( G \) graphs. Notably, a subfamily of \( G \) graphs exhibits a spectral redundancy index of \( 4/3 \), which matches that of the entire \( G\)-graph family. Other subfamilies of \( \mathcal{G} \) such as \( \mathcal{G}_{3q, q} \), \( \mathcal{G}_{4q, q} \), and beyond, though not analyzed in this study, also have finite spectral redundancy indices but these are strictly less than \( 4/3 \). While we do not examine these subfamilies in detail, we estimate that their indices are lower than \( 4/3 \). For instance, the spectral redundancy indices of \( \mathcal{G}_{3q, q} \) and \( \mathcal{G}_{4q, q} \) are \( 6/5 \) and \( 8/7 \) respectively.

Our study shows that the family of \( G \) graphs has a spectral redundancy index of \( 4/3 \). Alberto Seeger \cite{See} also examined the spectral redundancy index for various families of graphs. While some families were found to have infinite spectral redundancy indices, Seeger identified that the family of broom graphs has a finite spectral redundancy index of \( 4/3 \). This raises an interesting question for future research regarding the existence of a family of graphs with a finite spectral redundancy index greater than \( 4/3 \).\\

\noindent{\bf Acknowledgements.}  The research of Pawan Kumar is supported by CSIR, India as a SENIOR RESEARCH FELLOWSHIP, file No. 09/112(0669)/2020-EMR-I. The research of S. Pirzada is supported by National Board for Higher Mathematics (NBHM) research project number\\ NBHM/02011/20/2022.\\

\noindent{\bf Data availibility} Data sharing is not applicable to this article as no data sets were generated or analyzed during the current study.\\

\noindent{\bf Conflict of interest.} The authors declare that they have no conflict of interest.

{}
\end{document}